\documentclass[12pt,twoside]{article}

\usepackage{amsmath,amssymb}
\usepackage{amsthm,mathrsfs}
\usepackage{color,times}
\usepackage{hyperref}
\usepackage{textcomp}
\usepackage{sectsty}
\usepackage{verbatim}
 \def \no{\nonumber}
 \def \aint { \diagup\!\!\!\!\!\!\int}

 \allowdisplaybreaks
\begin{document}
 \footskip=0pt
 \footnotesep=2pt
\let\oldsection\section
\renewcommand\theequation{\arabic{equation}}
\newtheorem{claim}{\indent Claim}
\newtheorem{theorem}{\indent Theorem}
\newtheorem{lemma}{\indent Lemma}
\newtheorem{proposition}{\indent Proposition}
\newtheorem{definition}{\indent Definition}
\newtheorem{remark}{\indent Remark}
\newtheorem{corollary}{\indent Corollary}
\title{\Large \bf
 Remarks on the $\mathbf Q$-curvature flow}

\author{Xuezhang  Chen\thanks{ X. Chen: chenxuezhang@sina.com;},
Li Ma\thanks{L. Ma: nuslma@gmail.com;}
 and Xingwang  Xu \thanks{X. Xu: matxuxw@nus.edu.sg.}
\\
\small $^\ast$BICMR, Peking University, Beijing 100871, P. R. China
\\
\small $^\ast$Department of Mathematics \& IMS, Nanjing University,
Nanjing 210093, P. R. China
\\
\small $^\dag$Department of Mathematics, Henan Normal University,
Xinxiang 453007, P. R. China
\\
\small and
\\
\small $^\ddag$Department of Mathematics, National University of
Singapore, Block S17 (SOC1), \\ \small 10 Lower Kent Ridge Road,
Singapore 119076, Singapore}

\date{}
\maketitle
\begin{abstract}
{The main purpose of this short note is to point out that the
negative gradient flow for the prescribed $\mathbf Q$-curvature
problem on $S^n$ can be extended to handle the case that the
$\mathbf Q$-curvature candidate $f$ may change signs.}
\end{abstract}

\vskip .4in

\indent {\bf 1.} Various prescribing curvature problems on a
manifold can be restated as follows: given a smooth function $f$
defined on $M^n$ with the metric $g$, can one find a conformal metric
$g_u=e^{2u}g$ such that the aforementioned curvature is equal to
$f$? The typical example is the prescribing scalar curvature problem
with $(M^n, g) = (S^n, g_{S^n})$ where $g_{S^n}$ is the standard round
metric. In past several decades, this problem has attracted a lot of
attention. Recently, several groups of people are interested in the
prescribing $Q$-curvature problem. It is well known that
both problems are equivalent to solving certain partial differential
equations. When the background manifold is the standard sphere, the
non-compactness of the conformal group made the problem more
interesting to study. We refer readers to \cite{wx2} for more
background materials on this type of problem. Recall the prescribing
$Q$-curvature problem on $S^n$ is equivalent to the solvability of
the following equation
\begin{equation}\label{prescribed_Q-curvature}
P_n u+(n-1)!=f e^{nu} \text{~~on~~} S^n,
\end{equation}
where $P_n=P_{g_{S^n}}$ is $n$-th order Paneitz operator. Notice
that Equation (\ref{prescribed_Q-curvature}) has a variational
structure, hence the variational approach is a natural tool to
consider. Along this line, with many people's effort, several
sufficient conditions have been found to guarantee the existence of
solutions to (\ref{prescribed_Q-curvature}), for instance, see
\cite{cy}, \cite{wx1}-\cite{wx2} and references therein.

Recently, Brendle in \cite{bre2} introduced a flow method to study
the problem. It seems this new method is more promising. The first
and third authors of the current paper have adopted this method to
deal with the general higher order prescribing $Q$-curvature problem
on $S^n$ \cite{chxu1}. However, the positivity of the curvature
candidate $f$ plays an important role in their argument. We observe
that for prescribing Gaussian curvature problem on $S^2$, Hong and
Ma \cite{mahong} have verified that the positivity of the curvature
candidate $f$ can be removed. Some observation for the fourth-order
equation on $T^4$ with changing sign curvature candidate has also
appeared in \cite{gexu}. The purpose of this note is to point out
that, in fact, the positivity on $f$ is in general not necessary.
Before stating our main result, we give the definition of
non-degeneracy first. A smooth function $f$ defined on $S^n$ is
called non-degenerate if it satisfies
 $$(\Delta_{S^n}f)^2+|\nabla f|_{S^n}^2 \not = 0 \text{~~on~~} S^n.$$
Our main claim in this note can be stated as the following:
\begin{theorem}\label{main_Thm}
Let $n \geq 4$ be an even integer. Suppose $f: S^n \to \mathbb{R}$
is a sign changing smooth function with $\int_{S^n} f(x)
d\mu_{S^n}>0$. Assume in addition that $f$ admits only isolated
critical points with non-degeneracy in the set $\{x \in S^n;
f(x)>0\}$. Let
\begin{equation}\label{morse_ind}
\gamma_i=\sharp\{q \in S^n;f(q)>0, \;\; \nabla_{S^n} f(q)=0, \;\;
\Delta_{S^n} f(q) < 0, \;\; \text{ind}(f,q)=n-i\},
\end{equation}
where $\text{ind}(f,q)$ denotes the Morse index of $f$ at critical
point $q$. If the system of equations
\begin{equation}\label{Morse_eqn}
\gamma_0=1+k_0,\gamma_i=k_{i-1}+k_i,1 \leq i \leq n,k_n=0,
\end{equation}
has no non-negative integer solutions for $k_i$, then there exists a
solution to $Q$-curvature equation (\ref{prescribed_Q-curvature}).
\end{theorem}

Since the most part of the argument is the same as in \cite{chxu1},
here we only indicate how to overcome the difficulty arising from
the non-positivity of the curvature candidate $f$.

\vskip .2in \indent {\bf 2.} The first thing one needs to take care
of is the estimate of the normalized factor $\alpha(t)$. Before we
do this, let us set the stage first. Let $n=2m \geq 4$ be an even
integer and $\omega_n$ be the volume of the standard sphere $S^n$.
Let $f$ be a sign changing smooth function on $S^n$. Motivated by S.
Brendle \cite{bre2}, M. Struwe \cite{str} and Malchiodi-Struwe
\cite{mal_str}, the first and third authors of this note introduced
in \cite{chxu1} the following flow equation
\begin{equation} \label{flow}
2 u_t=\alpha f-Q,
\end{equation}
 where $Q=Q_g$ is the $Q$-curvature of the conformal metric
 $g(t)=e^{2u(t)}g_{S^n}$ which can be calculated by the formula
\begin{equation}\label{Q-curvature}
Qe^{nu}=P_n u+(n-1)!.
\end{equation}
Set
\begin{eqnarray*}
C_f^\infty=\left\{w \in C^\infty(S^n);g_w=e^{2w}g_{S^n}
\hbox{~~satisfies~~}\int_{S^n}d\mu_{g_w}=\omega_n
\hbox{~~and~~}\int_{S^n} f d\mu_{g_w}>0\right \}.
\end{eqnarray*}
We assume the flow (\ref{flow}) has the initial data $u(0) = u_0 \in
C_f^\infty$.

Recall, when $n$ is even, $P_n$ is given by
$$P_n=\prod_{k=0}^{(n-2)/2}(-\Delta_{S^n}+k(n-k-1)).$$

Observe that $P_n$ is a divergent free operator, hence integrating
(\ref{Q-curvature}) over $S^n$ yields
$$\aint_{S^n}Qe^{nu}d\mu_{S^n}=(n-1)!,$$
where $\aint_{S \sp n}$ denotes the average of the integral over
$S^n$. The energy functional associated with the equation
(\ref{flow}) can be written as
$$E_f[u]=E[u]-(n-1)!\log\left(\aint_{S^n}fe^{nu}d\mu_{S^n}\right),$$
where
$$E[u]={n \over 2}\aint_{S \sp n}(u P_nu+2(n-1)!u)d \mu_{S \sp n}.$$

We remark here that the flow (\ref{flow}) is the negative gradient
flow of the energy $E_f[u]$.

The normalized factor $\alpha(t)$ is chosen to be
\begin{equation}\label{alpha}
\alpha(t)=\frac{(n-1)!}{\aint_{S \sp n}f e^{nu}d \mu_{S^n}}.
\end{equation}
The reason to do so is to keep the volume of the flow metric $g(t)$
unchange for all time $t$, that is,
$$\aint_{S^n}e^{nu(t)}d\mu_{S^n} \equiv 1
\text{~~for all~~} t \geq 0.$$

In view of Lemma $1.1$ in \cite{chxu1}, the energy functional
$E_f[u(t)]$ is non-increasing, more explicitly, by a simple
calculation, one has
\begin{equation}\label{eng_dec}
\frac{d}{dt}E_f[u]=-{n \over 2}\aint_{S \sp
n}\left(\alpha(t)f-Q\right)\sp 2 d \mu_g.
\end{equation}

For sign changing curvature function $f$, similar to \cite{mahong},
we first need the following important observation.
\begin{lemma}\label{pres_C_f^infty}
If $u_0 \in C_f^\infty$, then for each time $t \geq 0$, the solution
$u(t)=u(t,u_0)$ is also in the class $C_f^\infty$. Moreover, there
exist two positive constants $C_1$ and $C_2$ depending only on $f$
and initial data $u_0$, such that
$$0 < C_1 \leq \alpha(t) \leq C_2 \hbox{~~for all~~} t \geq 0.$$
\end{lemma}
\noindent{\bf Proof.~~} By the selection of $\alpha(t)$, we first
need verify that if $u_0 \in C_f^\infty$, then $\int_{S^n}
fe^{nu(t)} d\mu_{S^n}>0$ for any time $t > 0$. In essence, with the
help of (\ref{eng_dec}) and Beckner's inequality (see \cite{beckner}
or \cite{chxu1} Prop. $1.1$), one has
\begin{eqnarray*}
-(n-1)!\log{(\max_{S^n}f)}& \leq & -(n-1)!\log{\aint_{S^n}f
e^{nu}d\mu_{S^n}}\\
&\leq& E_f[u](t)\leq E_f[u_0]<\infty.
\end{eqnarray*}
 Thus there hold
\begin{eqnarray*}
&&\max_{S^n} f \geq \aint_{S^n} f e^{n u(t)} d\mu_{S^n} \geq
e^{-E_f[u_0] \over (n-1)!}>0
\end{eqnarray*}
and
\begin{eqnarray*}
0 \leq E[u]& = & E_f[u]+(n-1)!\log\left(\aint_{S^n}fe^{nu}d\mu_{S^n}\right)\no\\
&\leq&E_f[u_0]+(n-1)! \log{(\max_{S^n}f)}<\infty.
\end{eqnarray*}
Furthermore, one can easily obtain:
\begin{equation}\label{bd_alpha}
{(n-1)! \over \max_{S^n}f} \leq \alpha(t) \leq (n-1)! e^{E_f[u_0]
\over (n-1)!}.
\end{equation}

Clearly Lemma \ref{pres_C_f^infty} follows from the equation
(\ref{bd_alpha}) with $C_1 = \frac{(n-1)!}{\max_{S^n} f}$ and $ C_2
= (n-1)! e^{E_f(u_0)/(n-1)!}$. \hfill $\Box$

With the help of this lemma, the integral estimate in \cite{chxu1}
go through without any changes.

\vskip .2in

\indent {\bf 3.} As usual, we have to investigate the property of
the compactness and the concentration behavior along the flow. To do
this, we follow the standard strategy to study its normalized flow
$v(t)$. It is well-known [for example, \cite{lyy}, Lemma
$5.4$~or~\cite{bre1}, Proposition $6$] that, for any family of
smooth functions $u(t)$, there exists a family of conformal
transformations $\phi(t): S^n \rightarrow S^n,$ smoothly depending
on the time $t$,
 such that
\begin{equation}\label{regular_constraint}
\aint_{S^n}x d\mu_h=0,
\text{ for all $t>0$,}
\end{equation}
with the normalized metric
\begin{equation}\label{normal_metric}
h=\phi^\ast (e^{2u} g_{S^n})\equiv e^{2v(t)}g_{S^n}.
\end{equation}

In view of the non-increasing property (\ref{eng_dec}) of
$E_f[u(t)]$ and a sharp version of Beckner's inequality (\cite{wx2},
Theorem $2.6$ or \cite{cy}), the global existence of the flow
(\ref{flow}) with any initial data $u_0 \in C_f^\infty$ is a direct
consequence of Section 2.1 in \cite{chxu1}.

We follow the same strategy as in the proof of Lemma 3.4 in
\cite{str} or Lemma 2.4 in \cite{chxu1} to obtain the asymptotic
behavior of $Q$-curvature of the flow metric, namely,
\begin{equation}\label{F_2_Q}
\int_{S^n}|\alpha f-Q|^2 d\mu_g \to 0 \hbox{~~as~~} t \to \infty.
\end{equation}

Then the rough curvature convergence (\ref{F_2_Q}) enables us to
employ Proposition $1.4$ of \cite{bre2} to the family of functions
$u_k=u(t_k)$ taking from the flow. We state it as the following:
\begin{lemma}\label{lem3.1}
Let $u_k=u(t_k), g_k=e^{2u_k}g_{S^n}$. Then, we have either (i) the
sequence $u_k$ is uniformly bounded in $H^n(S^n,g_{S^n})
\hookrightarrow L^\infty(S^n)$; or (ii) there exist a subsequence of
$u_k$ and finitely many points $q_1,\cdots ,q_L \in S^n$ such that
for any $r>0$ and any $l \in \{1,\cdots,L\}$, there holds
\begin{equation}\label{bre_con_points}
\liminf\limits_{k \to \infty}\int_{B_r(q_l)}|Q_k|d\mu_k\geq {1 \over
2}(n-1)! \omega_n,
\end{equation}
where $d\mu_k=d\mu_{g_k}$ and $Q_k=Q_{g_k}$ is the $Q$-curvature of
the metric $g_k$; in addition, the sequence $u_k$ is uniformly
bounded on any compact subset of $(S^n \setminus \{q_1,\cdots,
q_L\},g_{S^n})$ or $u_k \to -\infty$ locally uniformly away from
$q_1,\cdots ,q_L$ as $k \to \infty$.
\end{lemma}

\indent Just as some previous work has shown, a refined version of
Lemma $\ref{lem3.1}$ is much needed in late analysis.
\begin{lemma}\label{blow-up anal}
Let $u_k$ be the sequence of smooth functions on $S^n$ in Lemma
\ref{lem3.1}. In addition, there exists some sign changing smooth
function $Q_\infty$ defined on $S^n$, satisfying
$\|Q_k-Q_\infty\|_{L^2(S^n,g_k)}\to 0$ as $k \to \infty$. Let
$h_k=\phi_k^\ast(g_k)=e^{2v_k}g_{S^n}$ be the corresponding sequence
of normalized metrics given in (\ref{regular_constraint})-(\ref{normal_metric}). Then up
to a subsequence, either
\begin{enumerate}
\item[(i)]$u_k \to u_\infty$ in $H^n(S^n,g_{S^n})$ as $k \to \infty$, where
$g_\infty=e^{2u_\infty}g_{S^n}$ has Q-curvature $Q_\infty$, or
\item[(ii)] there exists $p \in S^n$, such that
\begin{equation}\label{concen_limit}
Q_\infty(p)=(n-1)! \hbox{~~and~~} d\mu_k \hookrightarrow \omega_n
\delta_p \text{~~as~~} k \to \infty,
\end{equation}
in the weak sense of measures, and
$$v_k \to 0 \text{~~in~~} H^n(S^n,g_{S^n}), Q_{h_k} \to (n-1)! \text{~~in~~} L^2(S^n,g_{S^n}).$$
In the latter case, $\phi_k$ converges weakly in
$H^{n/2}(S^n,g_{S^n})$ to the constant map $p$.
\end{enumerate}
\end{lemma}
\noindent{\bf Proof.~~} The proof follows the same argument as
Malchiodi and Struwe did in \cite{mal_str} or Chen and Xu in
\cite{chxu1}. When concentration occurs in the sense of
(\ref{bre_con_points}), we do need to overcome some difficulties
arising from the sign changing of $f$. For each $k$, choose $p_k \in
S^n$ and $r_k>0$ such that
\begin{equation}\label{r_k_p_k}
\sup\limits_{p \in
S^n}\int_{B_{r_k}(p)}|Q_k|d\mu_k=\int_{B_{r_k}(p_k)}|Q_k|d\mu_k={1
\over 4}(n-1)!\omega_n.
\end{equation}
Then by (\ref{bre_con_points}), $r_k \to 0$ as $k \to \infty$. Also
we may and will assume $p_k \to p$ as $k \to \infty$. For brevity,
one regards $p$ as $N$, the north pole of $S^n$.\\
\indent Denote by $\hat{\phi}_k$: $S^n \to S^n$ the conformal
diffeomorphisms mapping the upper hemisphere $S^n_+\equiv S^n \cap
\{x^{n+1}>0\}$ into $B_{r_k}(p_k)$ and the equatorial sphere
$\partial S^n_+$ to $\partial B_{r_k}(p_k)$. Indeed, up to a
rotation, $\hat{\phi}_k$ can be written as $\psi^{-p_k} \circ
\delta_{r_k} \circ \pi^{-p_k}$, where $\pi^{-p_k}: S^n \to
\mathbb{R}^n$ is the stereographic projection from $-p_k$ with the
inverse $\psi^{-p_k}=(\pi^{-p_k})^{-1}$ while the $\delta_{r_k}$ is
the dilation map on $\mathbb{R}^n$ defined by $\delta_{r_k}(y) =
\delta_{r_k} y$. In particular, set $\psi=\psi^{S}$. Consider the
sequence of functions $\hat{u}_k: S^n \to \mathbb{R}$ defined by
$$e^{2\hat{u}_k}g_{S^n}=\hat{\phi_k}^*(g_k)$$
which solve the equation
$$P_n \hat{u}_k+(n-1)!=\hat{Q}_k e^{n \hat{u}_k} \text{~~on~~} S^n,$$
where $\hat{Q}_k=Q_k \circ \hat{\phi}_k$. From the selection of
$r_k,p_k$ and (\ref{r_k_p_k}), by applying Lemma \ref{lem3.1} to
$\hat{u}_k$, we conclude that $\hat{u}_k \to \hat{u}_\infty$ in
$H^n_{\text{loc}}(S^n\setminus\{S\},g_{S^n})$ as $k \to \infty$,
where $S$ is the south pole on $S^n$. Meanwhile, $\hat{Q}_k \to
Q_\infty(p)$ almost everywhere as $k \to \infty$. Introducing the sequence of functions $\tilde{u}_k: S^n \to
\mathbb{R}$ by
$$e^{2\tilde{u}_k}g_{\mathbb{R}^n}=(\psi^{-p_k})^\ast(e^{2\hat{u}_k}g_{S^n})=\tilde{\psi}_k^\ast(g_k),$$
where $\tilde{\psi}_k=\psi^{-p_k} \circ \delta_{r_k}$, namely,
$$\tilde{u}_k=u_k \circ \tilde{\psi}_k+{1 \over n}\log (\det d \tilde{\psi}_k),$$
we find that $\tilde{u}_k$ converges in
$H^n_{\text{loc}}(\mathbb{R}^n)$ to a function $\tilde{u}_\infty$,
satisfying
\begin{equation}\label{eq_con_Q_cur}
(-\Delta_{\mathbb{R}^n})^{n/2}\tilde{u}_\infty=Q_\infty(p) e^{n
\tilde{u}_\infty}\text{~~in~~} \mathbb{R}^n.
\end{equation}
 Moreover, by Fatou's lemma we get
\begin{equation}\label{bd_vol}
\int_{\mathbb{R}^n}e^{n \tilde{u}_\infty}dz \leq \liminf\limits_{k
\to \infty} \int_{\mathbb{R}^n}e^{n\tilde{u}_k}dz=\omega_n.
\end{equation}
\indent Based on the proof of Lemma 3.2 of \cite{chxu1}, we need a
preliminary lemma to finish the proof of Lemma \ref{blow-up anal}.
\begin{lemma}\label{simple_bubble}
Under assumptions on $u_k$ as in Lemma \ref{blow-up anal}, there
holds $Q_\infty(p)>0$ and the solution $\tilde{u}_\infty$ to
equations (\ref{eq_con_Q_cur})-(\ref{bd_vol}) has the form
\begin{equation}\label{u_infty}
\tilde{u}_\infty(z)=\log{2\lambda \over 1+|\lambda (z-z_0)|^2}-{1
\over n}\log{Q_\infty(p) \over (n-1)!}
\end{equation}
for some $\lambda>0$ and $z_0 \in \mathbb{R}^n$.
\end{lemma}
\noindent{\bf Proof.~~} For brevity, one uses $u_\infty$ instead of
$\tilde{u}_\infty$. Let
$$\bar{w}(\rho)=\aint_{\partial B_\rho(0)}w(z) d\sigma(z), \rho>0$$
denote the spherical average of the function $w$ defined in
$\mathbb{R}^n$. Due to the proof of Lemma 3.3 in \cite{chxu1}, we
only need rule out the case of $Q_\infty(p)\leq 0$. Arguing by
negation, we assume $Q_\infty(p)\leq 0$. The argument heavily relies
on the following important estimate obtained through the analysis on
Green's function of some $(n-2)$-order differential operator in
\cite{chxu1}, Lemma $3.3$. For convenience, we restate it here: for
any $r>0$ and $q \in S^n$, there holds
\begin{equation}\label{est1}
|\int_{B_r(q)}\Delta_{S^n} u_k d\mu_{S^n}|\leq B_0 r^{n-2}
\end{equation}
for all $k$, where $B_0>0$ is a constant.

\indent Let $m={n/2} \geq 2$ and $w_i(z)=(-\Delta)^{i}u_\infty(z),
i=1,2,\cdots,m$. Then, we claim that for $1 \leq i \leq m-1$, there
holds
\begin{equation}\label{ineq6}
w_{m-i} \leq 0 \text{~~in~~} \mathbb{R}^n.
\end{equation}
For $w_{m-1}$, by negation, we assume there exists $z_0 \in
\mathbb{R}^n$, such that $w_{m-1}(z_0)>0$. Without loss of
generality, assume $z_0=0$. From (\ref{eq_con_Q_cur}) and Jensen's
inequality, we have
\begin{eqnarray}
 \left \{
\begin {array}{llll}
&-\Delta \bar{u}_\infty=\bar{w}_1,\\
&-\Delta \bar{w}_1=\bar{w}_2,\\
&~~~~~\cdots \\
&-\Delta \bar{w}_{m-1}=\bar{w}_m \leq Q_\infty(p) e^{n
\bar{u}_\infty}\leq 0.
\end{array}
\right.\label{system}
\end{eqnarray}
Thus $\bar{w}_{m-1}'(\rho)\geq 0$, which indicates
$\bar{w}_{m-1}(\rho)\geq \bar{w}_{m-1}(0)=w_{m-1}(0)>0$. Observe
that
$$\bar{w}_{m-2}'(\rho)={-\rho \over n}[|B_\rho(0)|^{-1} \int_{B_{\rho}(0)}w_{m-1}(z)dz] \leq {(-w_{m-1}(0))\over n}\rho<0.$$
Thus it follows that
$$-\bar{w}_{m-2}(\rho)\geq B_2 \rho^2 \text{~~for~~} \rho \geq \rho_1>0, B_2>0.$$
By (\ref{system}) and mathematical induction, in general, for $2 \leq i \leq m-1$, we have
$$(-1)^{i-1} \bar{w}_{m-i}(\rho)\geq B_i \rho^{2(i-1)} \text{~~for~~} \rho \geq \rho_{i-1}>0, B_i>0.$$
Apply this to $i=m-1$ to get
\begin{equation}\label{ineq9}
(-1)^{m-2} \int_{\partial B_\rho(0)}(-\Delta u_\infty(z)) d\sigma(z)
\geq B_{m-1}\rho^{2(m-2)+n-1} \text{~~for~~} \rho \geq \rho_{m-1}.
\end{equation}
For sufficiently large $k$ and all $\rho \geq \rho_{m-1}$, one has
\begin{equation}\label{ineq17}
(-1)^{m-2}\int_{\partial B_\rho(0)}(-\Delta \tilde{u}_k(z))dz \geq
A_1 \rho^{2(m-2)+n-1}
\end{equation}
where $A_1>0$ is a universal constant. By a similar argument on
pages 951-953 of \cite{chxu1}, using (\ref{ineq17}) and the
expression
$$\tilde{u}_k(z)=u_k \circ \tilde{\psi}_k+\log{2r_k \over 1+|r_k z|^2},$$
one obtains that for some fixed $L>0$ and any $d>L$, there holds
\begin{equation}\label{ineq7}
(-1)^{m-2}\int_{B_{dr_k}(p_k)}(-\Delta_{S^n}u_k )d\mu_{S^n} \geq A_2
r_k^{n-2}(d^{2(m-2)+n}-L^{2(m-2)+n}-L^{n-2})
\end{equation}
for sufficiently large $k$,where $A_2>0$ is a constant.\\
\indent On the other hand, by choosing $r=r_k d$ and $q=p_k$ in
(\ref{est1}), with a uniform constant $A_3>0$ it yields
\begin{equation}\label{ineq8}
(-1)^{m-2}\int_{B_{d r_k}(p_k)}(-\Delta_{S^n} u_k) d\mu_{S^n} \leq
A_3 r_k^{n-2} d^{n-2}.
\end{equation}
Hence, for any fixed $L>0$ as above and sufficiently large $k$, (\ref{ineq7}) and (\ref{ineq8})
yield a contradiction by choosing $d$ sufficiently large. \\
\indent Next, we prove (\ref{ineq6}) by the induction argument. The
case $i=1$ has been settled above. If $m = 2$, we are done. Thus we
assume $m
> 2$. Suppose for some $i$ with $1 \leq i \leq m-2$ and all $1 \leq
k \leq i$, $w_{m-k} \leq 0$ in $\mathbb{R}^n$. Then one needs to
show $w_{m-i-1}\leq 0$ in $\mathbb{R}^n$. By negation again, we may
assume $w_{m-i-1}(0)>0$. Since $\bar{w}_{m-i-1}'(\rho)={-1 \over
|\partial B_\rho(0)|}\int_{B_\rho(0)}w_{m-i}(z)dz \geq 0$, it
follows that
$$\bar{w}_{m-i-1}(\rho) \geq \bar{w}_{m-i-1}(0)=w_{m-i-1}(0)>0.$$

If $ i \le m - 3$, by $-\Delta \bar{w}_{m-i-2}=\bar{w}_{m-i-1}$, one
has
$$-\bar{w}_{m-i-2}(\rho) \geq B_2 \rho^2 \text{~~for~~} \rho \geq \rho_1>0, B_2>0.$$
In general, by (\ref{system}) one obtains
$$(-1)^{j+1} \bar{w}_{m-i-j}(\rho)\geq B_j \rho^{2(j-1)} \text{~~for~~} \rho \geq \rho_{j-1}>0,B_j>0,i+j\leq m-1.$$
\indent  Choosing $j=m-1-i$, we have
\begin{equation}\label{ineq10}
(-1)^{m-i}\int_{\partial B_\rho(0)}(-\Delta u_\infty(z))d\sigma(z)
\geq B_{m-1-i}\rho^{2(m-i-2)+n-1} \text{~~for~~}\rho \geq
\rho_{m-2-i}.
\end{equation}
Fixing $L\geq \rho_{m-2-i}$ and for any $d>L$, by a similar argument
on (\ref{ineq7}), with a constant $A_5>0$, one gets
\begin{equation}\label{ineq11}
(-1)^{m-i}\int_{B_{dr_k}(p_k)}(-\Delta_{S^n}u_k) d\mu_{S^n}\geq A_5
r_k^{n-2} (d^{2(m-i-2)+n}-L^{2(m-i-2)+n}-L^{n-2})
\end{equation}
for all sufficiently large $k$. On the other hand, choosing $r=r_k
d$ and $q=p_k$ in (\ref{est1}), with another constant $A_6>0$ one
has
\begin{equation}\label{ineq12}
(-1)^{m-i}\int_{B_{dr_k}(p_k)}(-\Delta_{S^n}u_k) d\mu_{S^n} \leq A_6
r_k^{n-2} d^{n-2}.
\end{equation}
Thus (\ref{ineq11}) and (\ref{ineq12}) would contradict each other if $d$ is chosen to be sufficiently large.\\
\indent If $i=m-2$, by inductive assumption, $w_2(z) \leq 0$ in
$\mathbb{R}^n$ and $w_1(0)>0$. Since $\bar{w}_1'(\rho) \geq 0$,
$$\bar{w}_1(\rho)\geq \bar{w}_1(0)=w_1(0)>0.$$
Given any $d>0$, from the above inequality, it is easy to derive
$$\int_{\partial B_\rho (0)}(-\Delta u_\infty) dz \geq 2 A_7 \rho^{n-1} \text{~~for~~} 0 \leq \rho \leq d,$$
where $A_7>0$ depends on $w_1(0)$ and $n$ only. Repeating the
argument for (\ref{ineq7}) and scaling back to $S^n$, one gets
\begin{equation}\label{ineq13}
\int_{B_{dr_k}(p_k)}(-\Delta_{S^n}u_k) d\mu_{S^n}\geq A_7 r_k^{n-2}
d^n,
\end{equation}
for all sufficiently large $k$. Now, the equation (\ref{est1}) with
$r= r_k d$ and $q=p_k$ gives
\begin{equation}\label{ineq14}
\int_{B_{dr_k}(p_k)}(-\Delta_{S^n}u_k) d\mu_{S^n} \leq A_8 r_k^{n-2}
d^{n-2}
\end{equation}
where $A_8>0$ is a constant. Equations (\ref{ineq13}) and (\ref{ineq14}) contradict each other if $d>0$ is sufficiently large.\\
\indent Therefore, we conclude that $w_{m-i-1} \geq 0$ in
$\mathbb{R}^n$ and the induction argument is complete.\\
\indent Finally, it follows from the inequality (\ref{ineq6}) that
$-\Delta u_\infty = w_1 \leq 0$ in $\mathbb{R}^n$, that is,
$u_\infty$ is a subharmonic function. By the mean value property for
subharmonic functions, for any $z \in \mathbb{R}^n$ and $r>0$, there
holds
$$n u_\infty(z)\leq  |B_r(0)|^{-1}\int_{B_r(z)}nu_\infty(y)dy.$$
By this inequality and Jensen's inequality, one gets
\begin{eqnarray}
e^{nu_\infty(z)} &\leq&
|B_r(0)|^{-1}\int_{B_r(z)}e^{nu_\infty(y)}dy\no\\
&\leq& |B_r(0)|^{-1}
\int_{\mathbb{R}^n}e^{nu_\infty(y)}dy\label{ineq18}
\end{eqnarray}
for any $r>0$. In view of  (\ref{bd_vol}), by letting $r \to
\infty$, the inequality (\ref{ineq18}) indicates that $u_\infty
\equiv -\infty$ in $\mathbb{R}^n$, which is obviously impossible.
The proof is complete. \hfill $\Box$

\vskip .2in

\noindent {\bf Proof of Lemma \ref{blow-up anal} (completed).} With
the help of Lemma \ref{simple_bubble}, we can show that there is the
unique concentration point $p$ of $\{g_k\}$ such that
$Q_\infty(p)=(n-1)!$. To see this, first set
$Q_\infty^+=\max{\{Q_\infty, 0\}}, Q_\infty^-=\min{\{Q_\infty,0\}}$.
Notice that by Lemma \ref{lem3.1}, there are only finitely many
blow-up points, say $p_1, p_2, \cdots, p_l$. By previous two Lemmas,
we know that at each $p_i$, $Q_\infty(p_i) > 0$. Now for each $i$,
choose a sufficiently small $r_i > 0$ so that $Q_\infty \geq 0$ in
$B_{r_i}(p_i)$ and $B_{r_i}(p_i)\cap B_{r_j}(p_j) = \emptyset$ if $
i \not = j$. Then follow the same argument on page 957 of
\cite{chxu1} (or similar one in \cite{mal_str}) to conclude that
\begin{equation}\label{cenpt1}
(n-1)!l \omega_n=\sum_{i=1}^l \int_{\mathbb{R}^n} Q_\infty(p_i)
e^{n\tilde{u}_\infty}dz \leq \sum_{i=1}^l \int_{B_{r_i}(p_i)}
Q_\infty^+ d\mu_k+o(1),
\end{equation}
for all sufficiently large $k$, where $o(1) \to 0$ as $k \to
\infty$. Thus, there holds $\lim_{k \to \infty}\int_{S^n}Q_\infty^-
d\mu_k=\sum_{i=1}^lQ_\infty^-(p_i)\omega_n=0$ since concentration
phenomena only occur at points $p_i$ where $Q_\infty(p_i)>0, 1 \leq
i \leq l$. From this identity and the selection of $r_i$, one has
\begin{eqnarray}
\sum_{i=1}^l \int_{B_{r_i}(p_i)} Q_\infty^+ d\mu_k & = &\sum_{i=1}^l
\int_{B_{r_i}(p_i)}Q_\infty
d\mu_k\no\\
& = & \sum_{i=1}^l [\int_{B_{r_i}(p_i)} Q_k d\mu_k+ \int_{B_{r_i}(p_i)}(Q_\infty -Q_k) d\mu_k]\no\\
& \leq & \int_{S^n} Q_k d\mu_k + 2\int_{S^n} |Q_\infty - Q_k| d\mu_k
+ \int_{S^n\backslash \cup_{i=1}^l B_{r_i}(p_i)} |Q_\infty| d\mu_k\no\\
&=&(n-1)!\omega_n+o(1),\label{cenpt2}
\end{eqnarray}
for all sufficiently large $k$ and where we have used the local
volume concentration property in the last term and uniform bound of
$Q_\infty$. Thus, it follows from \eqref{cenpt1} and \eqref{cenpt2}
that $l = 1$ and $Q_\infty(p)=(n-1)!$. Finally, the rest part of the
proof of Lemma \ref{blow-up anal} is the same as the proof of Lemma
3.2 in \cite{chxu1}. \hfill $\Box$
\begin{remark}
We should point out that, one can not apply Theorem 9 in \cite{mar}
to derive Lemma \ref{blow-up anal} directly. The assumption in
\cite{mar}: $Q_k \to Q_\infty$ in $C^0(S^n)$ is much stronger than
the one in Lemma \ref{blow-up anal}. Similar blow-up analysis as in
\cite{mar} has also been done by Malchiodi \cite{mal}. However those
estimates seem not suitable for Q-curvature flow since it is hard to
have $C^0$ convergence. So we have to seek another reasonable
procedure to do blow-up analysis in the flow setting.
\end{remark}
The remainder of the proof of Theorem \ref{main_Thm} will be
completed through a contradictive argument. From now on, we assume
$f$ can not be realized as a $Q$-curvature of any metric in the
conformal class of $g_{S^n}$. Following the standard scheme in
\cite{chxu1}, in particular Sections 4-5, along with Lemma
\ref{blow-up anal}, one eventually obtains the asymptotic behavior
of the flow $u(t)$ and the so-called shadow flow
$$\Theta = \Theta(t)= \aint_{S^n} \phi(t) d\mu_{S^n}.$$
\begin{lemma}\label{limit_E_f[u]}
Let $u(0)=u_0 \in C_f^\infty$ be the initial data of the flow
(\ref{flow}) and (\ref{Q-curvature}). Then the flow metrics $g(t)$
concentrate at a critical point $p$ of $f$ with $f(p)>0,
\Delta_{S^n}f(p)\leq 0$ and the energy $E_f[u(t)]$ converges to
$-(n-1)! \log f(p)$, that is
$$E_f[u(t)] \to -(n-1)! \log f(p), \text{ as } t \to \infty.$$
Moreover,  the critical point $p$ is also the unique limit of the
shadow flow $\Theta(t)$ associated with $u(t)$, in other words, $ p
= \lim_{t \to \infty}\Theta(t)$.
 \end{lemma}

\vskip .2in

\indent {\bf 4.} In this and next part, we will briefly sketch the
proof of our main result. For $q \in S^n, 0 < \epsilon <\infty,$
denote by $\phi_{-q,\epsilon}=\psi^{-q}\circ \delta_\epsilon \circ
\pi^{-q}$ the stereographic projection with $-q$ at infinity, that
is, $q$ becomes the north pole in the stereographic coordinates. It
is relatively easy to verify that $\phi_{-q,\epsilon}$ converges
weakly in $H^{n/2}(S^n,g_{S^n})$ to $q$ as $\epsilon \to 0$. Define
a map
$$j: S^n\times (0,\infty)\ni (q,\epsilon) \mapsto u_{q,\epsilon}=
{1 \over n}\log \det(d\phi_{q,\epsilon})\in C_\ast^\infty.$$
 And set
$g_{q,\epsilon}=\phi_{q,\epsilon}^\ast(g_{S^n})=e^{2u_{q,\epsilon}}g_{S^n}.$
Then we have
$$d\mu_{g_{q,\epsilon}}=e^{n u_{q,\epsilon}}d\mu_{S^n} \rightharpoonup \omega_n \delta_q,$$
in the weak sense of measures as $\epsilon \to 0$. For $\gamma \in
\mathbb{R},$ denote by
$$L_{\gamma}=\{u\in C_f^\infty; E_f[u] \leq \gamma\},$$
the sub-level set of $E_f$. Under our assumptions on $f$, label all
critical points of $f$ with positive critical values by $q_1,
\cdots,q_N$ such that $0<f(q_i) \leq f(q_j)$ for $1 \leq i \leq j
\leq N$ and set
$$\beta_i=-(n-1)!\log f(q_i)=\lim\limits_{\epsilon \to 0}E_f[u_{q_i,\epsilon}], 1 \leq i \leq N.$$
Without loss of generality, we assume all critical levels $f(q_i),1
\leq i \leq N$ are distinct, so there exists a $\nu_0>0$ such that
$\beta_i-2\nu_0 > \beta_{i+1}$, in fact we can take $\nu_0 =
\frac{1}{2} \min_{1 \le i \le N - 1} \{\beta_i - \beta_{i+1}\} > 0$.\\
\indent First of all, we shall characterize the homotopy types of
the sub-level sets. We state them as a proposition, which has
analogous counterpart in \cite{mal_str} or \cite{chxu1}.

\begin{proposition}\label{morse_prop}
\begin{enumerate}
\item[(i)] If $\delta_0 > \hbox{max}\{-(n-1)! \log (\aint_{S^n}f(x)d\mu_{S^n}),\beta_1\} $, the set $L_{\delta_0}$ is contractible.

\item[(ii)] For any $0 < \nu \leq \nu_0$ and each $1 \leq i \leq N$, the sets
$L_{\beta_i-\nu}$ and $L_{\beta_{i+1}+\nu}$ are homotopy equivalent.

\item[(iii)] For each critical point $q_i$ of $f$ where
$\Delta_{S^n}f(q_i)>0$ and $f(q_i)>0$, the sets $L_{\beta_i+\nu_0}$
and $L_{\beta_i-\nu_0}$ are homotopy equivalent.

\item[(iv)] For each critical point $q_i$ where $\Delta_{S^n}f(q_i)<0$ and $f(q_i)>0$, the
set $L_{\beta_i+\nu_0}$ is homotopic to the set $L_{\beta_i-\nu_0}$
with $(n-\text{ind}(f,q_i))$-cell attached.
\end{enumerate}
\end{proposition}
\noindent{\bf Proof:} (i) Let $\delta_0$ be chosen as above. For $0
\leq s \leq 1$ and $u_0 \in C_f^\infty$, define
$$H_1(s,u_0)={1 \over n}\log((1-s)e^{n u_0}+s), \text{~~that is~~} e^{n H_1(s,u_0)}=(1-s)e^{nu_0}+s,$$
then one easily obtains
\begin{eqnarray*}
&&\aint_{S^n} e^{n H_1(s,u_0)}d\mu_{S^n}=1 \hbox{~~and~~}\\
&&\int_{S^n}f e^{n H_1(s,u_0)}d\mu_{S^n}=(1-s)\int_{S^n} f
e^{nu_0}d\mu_{S^n}+s \int_{S^n}f d\mu_{S^n}>0,
\end{eqnarray*}
in view of the assumption that $\int_{S^n}f(x)d\mu_{S^n}>0$,
$H_1(s,u_0)$ provides a homotopic deformation within the set
$C_f^\infty$. Given such $u_0$ and $0 \leq s \leq 1$, by Lemma
$\ref{limit_E_f[u]}$ and the selection of $\delta_0$, there exists a
minimal time $T=T(s,u_0)$, such that $E_f[u(T,H_1(s,u_0))] \leq
\delta_0$, where the continuity of $T(s,u_0)$ on $s$ and $u_0$ can
be deduced by (\ref{eng_dec}) and the expression of $H_1(s,u_0)$.
Thus the map
 $H:(s,u_0) \mapsto u(T(s,u_0),H_1(s,u_0))$
is the desired contraction of $L_{\delta_0}$ within itself.  To see
this, first, by lemma \ref{pres_C_f^infty}, one knows that $H(s,u_0)
\in C_f^\infty$; next notice that $T(0,u_0) = 0$, hence $u(T(0,
u_0), H(0, u_0)) = u(0, u_0) = u_0$ and
$u(T(1, u_0), H(1, u_0)) = 0$ since $H(1, u_0) = 0$ with $E_f[0]=-(n-1)! \log (\aint_{S^n}f(x)d\mu_{S^n}) < \delta_0, T(1, u_0) = 0$.\\
\indent The proofs of (ii)-(iv) are identical to the corresponding
ones of Proposition $6.1$ (ii)-(iv) in \cite{chxu1}.\hfill $\Box$

\vskip .2in

\indent {\bf 5.} We are now in position to complete the proof of our
main theorem.

\vskip .2in

\noindent{\bf Proof of Theorem \ref{main_Thm}:} By negation, suppose
the flow is divergent for any initial data in $C_f^\infty$ and there
is no conformal metric of $g_{S^n}$ with $Q$-curvature $f$.
Proposition \ref{morse_prop} shows that for some suitable
$\delta_0$, $L_{\delta_0}$ is contractible and homotopically
equivalent to a set $E_\infty$ whose homotopy type consists of a
point $\{p\}$ with $(n-\text{ind}(f,q))$-dimensional cell attached
for each critical point $q$ of $f$ with $\Delta_{S^n}f(q)<0$ and
$f(q)>0$. By applying \cite{changkc}, Theorem $4.3$ on page $36$ to
$L_{\delta_0}$, we conclude that the identity
\begin{equation}\label{morse_chang}
\sum\limits_{j=0}^n s^j \gamma_j=1+(1+s)\sum\limits_{j=0}^n s^j k_j
\end{equation}
holds for Morse polynomials of $L_{\delta_0}$ and $E_\infty$, where
$k_j \ge 0$ are integers and $\gamma_j$ is defined in
(\ref{morse_ind}). Thus we achieve a contradiction with the
assumption that the system (\ref{Morse_eqn}) has no nonnegative
integer solutions $k_j$ and this contradiction completes the proof.
\hfill $\Box$

\vskip .2in

\noindent {\bf Acknowledgments:} We would like to thank the referees
for critical comments. The first author is partially supported
through NSF of China (No.11201223) and China postdoctoral foundation
(No.2011M500175). First author would like to thank Math Department
of NUS for their hospitality and financial support, and is grateful
to Professor Michael Struwe for stimulating discussions by emails
and bringing reference \cite{mar} to his attention. The second
author's research is partially supported by the National Natural
Science Foundation of China (No.11271111). The third author's
research is partially supported by NUS research grant
R-146-000-127-112 as well as the Siyuan foundation through Nanjing
University.

\noindent {\bf Added notes by Xuezhang Chen:} This is a second version of this paper, which is essentially completed in January of 2013. The first version is completed in April of 2012.

\end{document}